\RequirePackage{fix-cm}
\documentclass[smallextended]{svjour3}       
\smartqed  
\usepackage{graphicx}
\usepackage{amssymb}
\usepackage{amsmath}
\usepackage{bbding}
\begin{document}

\title{The weighted horizontal linear complementarity problem on a Euclidean Jordan algebra
}

\titlerunning{weighted complementarity problem}        

\author{Xiaoni Chi         \and
        M. Seetharama Gowda \and
        Jiyuan Tao 
}


\institute{Xiaoni Chi \at
School of Mathematics and Computing Science,
Guangxi Colleges and Universities Key Laboratory of Data Analysis and Computation,
Guilin University of Electronic Technology,
Guilin {\rm 541004},  Guangxi,  P. R. China \\
\email{chixiaoni@126.com}           
           \and
           M. Seetharama Gowda (\Envelope)\at
Department of Mathematics and Statistics,
University of Maryland, Baltimore County,
Baltimore, Maryland  \rm{21250}, USA\\
\email{gowda@umbc.edu}
           \and
           Jiyuan Tao \at
Department of Mathematics and Statistics,
Loyola University Maryland,
Baltimore, Maryland \rm{21210}, USA\\
\email{jtao@loyola.edu}
}

\date{Received: date / Accepted: date}

\maketitle

\begin{abstract}
A weighted complementarity problem (wCP) is to find a pair of vectors belonging to the
intersection of a manifold and a cone such that the product of the vectors in a certain  algebra equals
a given weight vector.
If the weight vector is zero, we get a complementarity problem.
Examples of such problems include
the Fisher market equilibrium problem and  the linear programming and weighted centering problem.
In this paper we consider the weighted horizontal linear complementarity problem (wHLCP) in the
setting of  Euclidean Jordan algebras and establish some existence and uniqueness results.
For a pair of linear transformations on a Euclidean Jordan algebra, we introduce the
concepts of ${\bf R}_0$, ${\bf R}$, and ${\bf P}$ properties and discuss the solvability of
wHCLPs under  nonzero (topological) degree conditions.
A uniqueness result is stated in the setting of $\mathbb{R}^{n}$. We show how our
results naturally lead to interior point systems.

\end{abstract}

\section{Introduction}
\label{intro}
Introduced in \cite{potra1}, a {\it weighted complementarity problem} (wCP)  is to find a pair of
vectors $(x,y)$ belonging to the intersection of a manifold with a cone such that their product in a
certain (Euclidean Jordan) algebra equals a given weight vector $w$. When $w$ is the zero vector, wCP becomes a
complementarity problem (CP).
To elaborate, consider a Euclidean Jordan algebra  $(V,\circ, \langle \cdot, \cdot\rangle)$  with
symmetric cone $V_+$ \cite{faraut-koranyi}. Given a map $F:V\times V\times \mathbb{R}^{l}\rightarrow V\times \mathbb{R}^{l}$ and 
a weight vector $w\in V_+$, wCP is to find $(x,y)\in V\times V$  such that for some $u$,
\begin{eqnarray}\label{wCP}
\begin{array}{l}
 x\in V_+,\ y\in V_+,\\
 x\circ y=w,\\
F(x,y,u)=0.
 \end{array}
\end{eqnarray}
Deferring this general problem for a future study, in \cite{potra1} and \cite{potra2}, Potra studies affine wCP on the (Euclidean Jordan) algebra
$\mathbb{R}^{n}$ with several examples and  results.
Given matrices $A,B\in \mathbb{R}^{(n+m)\times n}$, $C\in \mathbb{R}^{(n+m)\times m}$, a weight vector $w\in\mathbb{R}^{n}_{+}$, and $q\in \mathbb{R}^{n+m}$,
the  {\it weighted mixed horizontal linear complementarity problem}  considered in \cite{potra1}, \cite{potra2} is to find  $(x,y,z)\in \mathbb{R}^{n}\times \mathbb{R}^{n}\times \mathbb{R}^{m}$ such that

\begin{eqnarray}\label{affine wCP}
\begin{array}{l}
 x\geq0,\ y\geq0,\\
 x\ast y=w,\\
 Ax+By+Cz=q,
 \end{array}
\end{eqnarray}
where $x*y$ denotes the Hadamard (= componentwise) product of vectors $x$ and $y$. Here, $C$ is assumed
to be of full column rank and so, with a suitable change of variables, one could transform
(see Section 2 in \cite{potra2}) the above
problem to an equivalent affine wCP where $C$ becomes vacuous and $A$ and $B$ are square.
It is shown in \cite{potra1} that
the Fisher  market equilibrium problem \cite{eisenberg-gale},\cite{ye} and
the linear  programming and weighted centering  problem \cite{anstricher}
can be formulated in the form (\ref{affine wCP}), where the triple $(A,B,C)$ satisfies a certain monotonicity condition.
In \cite{potra1}, Potra presented and analyzed two interior-point methods for
solving such  monotone affine wCPs.
Subsequently, replacing `monotone' conditions by `row and column sufficient' conditions,
Potra \cite{potra2}
described several theoretical results and  a corrector-predictor interior-point method for its numerical solution.

We note that weighted complementarity problems  were studied much earlier in connection with interior point methods.
For example, in \cite{kojima et al}, Kojima et al showed that if (continuous) $f:\mathbb{R}^{n}\rightarrow \mathbb{R}^{n}$ is a uniform P-function, that is, there exists some $\gamma>0$ such that
$$\max_{1\leq i\leq n} (x_i-y_i)\Big ( f_i(x)-f_i(y)\Big )\geq \gamma ||x-y||^2 \,\,\,(\mbox{for all}\,\,x,y\in \mathbb{R}^{n}_+),$$
then the mapping
$\widehat{f}: \mathbb{R}^{n}_{+}\times \mathbb{R}^{n}_{+}\rightarrow \mathbb{R}^{n}_{+}\times \mathbb{R}^{n}$ is a homeomorphism, where
$$
\widehat{f}(x,y)=\left [ \begin{array}{c}
x*y\\y-f(x)\end{array}\right ].$$
 In this situation, the following
weighted nonlinear complementarity problem has a unique solution for each $w\in \mathbb{R}^{n}_{+}$ and
$q\in \mathbb{R}^{n}$:

\begin{eqnarray}
\begin{array}{l}
 x\geq0,\ y\geq0,\\
 x\ast y=w,\\
 y=f(x)+q.
 \end{array}
\label{f3}
\end{eqnarray}

Another work that specifically looks at the general wCP (\ref{wCP}) is by Yoshise \cite{yoshise}. In this work,  under certain`monotonicity and injectivity' assumptions, it is shown that the map
$(x,y,u)\rightarrow \Big (x\circ y, F(x,y,u)\Big )$ is a homeomorphism on a certain subset of $V_+\times V_+\times \mathbb{R}^{l}$, leading to the solvability of wCP, see Theorem 3.10 and Corollary 4.4 in \cite{yoshise}.

Our objective in this paper is to consider the following affine wCP in the setting of Euclidean Jordan algebras. Given
two linear transformations
$A$ and $B$ on a Euclidean Jordan algebra $V$, a weight vector $w\in V_+$, and $q\in V$, the
{\it weighted horizontal linear complementarity problem} wHLCP$(A,B,w,q)$
is to find
$(x,y)\in V\times V$  such that
\begin{eqnarray}\label{whlcp}
\begin{array}{l}
 x\geq0,\ y\geq0,\\
 x\circ y=w,\\
 Ax+By=q,
 \end{array}
\label{f1}
\end{eqnarray}
where $x\geq 0$ means that $x\in V_+$, etc.
If $w=0$, the above problem reduces to the (symmetric cone) {\it horizontal linear complementarity problem} on $V$, denoted by  HLCP$(A,B,q)$.
If $w=0$, $A=I$, and $B=-M$, HLCP$(A,B,q)$ reduces to the  (symmetric cone) {\it linear complementarity problem} LCP$(M,q)$ on $V$.
In particular, when $V=\mathbb{R}^{n}$, this reduces to the {\it standard linear complementarity problem}.

Our analysis differs from Potra's \cite{potra1}, \cite{potra2} in several ways. First, our setting is that of a
general Euclidean Jordan algebra, instead of $\mathbb{R}^{n}$. Second, instead of the `monotone/sufficient' conditions, we rely on the ${\bf R}_0$ property (that is commonly
used
 in the LCP literature) coupled with a nonzero degree condition of a certain map associated with wHLCP.
Third, instead of using the optimization methodology, we rely on the degree theoretic tools.
Our analysis also differs from that of Yoshise \cite{yoshise} where results were proved under certain
`monotonicity and injectivity' conditions.

In this paper, we establish some basic  existence/uniqueness results about
wHLCPs.
Generalizing the LCP concept of a degree of an ${\bf R}_0$-matrix, we introduce the concept of
degree of an ${\bf R}_0$-pair of linear transformations in the setting of Euclidean Jordan algebras.
Assuming that this degree is nonzero for the pair $\{A,B\}$, we show that  wHLCP$(A,B,w,q)$ has a nonempty compact solution set
for every  $(w,q)\in V_+\times V.$
This conclusion, in particular, will allow us to say that
\begin{itemize}
\item [$\bullet$] the map $\Gamma:V_+\times V_+\rightarrow V_+\times V$ given by
$\Gamma(x,y)=\Big ( x\circ y, Ax+By\Big )$ is surjective and
\item [$\bullet$] when $w>0$ (that is, $w\in \mbox{int}(V_+)$),  the `interior point system'
$$x>0,\,y>0,\,\,x\circ y=w,\,\,\mbox{and}\,\,Ax+By=q$$ has a nonempty compact solution set.
\end{itemize}
The result about the interior point system appears to be new even in the setting of standard LCPs.

We also introduce the concept of a  ${\bf P}$-pair and show that
when $V=\mathbb{R}^{n}$, wHLCP$(A,B,w,q)$ has a unique solution
for every  $(w,q)\in\mathbb{R}^{n}_+\times\mathbb{R}^{n}$.

The organization of our paper is as follows. In Section 2, we cover some
basic material. In Section 3, we introduce the concepts of ${\bf R}_0$ and ${\bf R}$ pairs and define the degree of an ${\bf R}_0$-pair.
Section 4 covers the main result of the paper describing the solvability of
wHLCP. While Section 5 deals with the solution set behavior, Sections 6 and 7 cover ${\bf P}$-pairs and address uniqueness issues.

\section{Preliminaries}
\label{sec:1}
Throughout this paper, $\mathbb{R}^{n}$
 denotes the Euclidean $n$-space of  real column
vectors. We use the (same) symbol $0$ to denote the  zero vector in any vector space.
$(V,\circ, \langle \cdot,\cdot\rangle)$ denotes a Euclidean Jordan algebra of rank $n$ with symmetric cone $V_+$ \cite{faraut-koranyi}, \cite{gowda-sznajder-tao}. Here, $x\circ y$ and $\langle x,y\rangle$, respectively, denote the Jordan product and the inner product of elements $x$ and $y$. The unit element of $V$ is denoted by $e$.  For a subset $S$ of $V$, the interior, closure, and boundary are denoted by
$int(S)$, $\overline{S}$, and $\partial(S)$. If $x\in V_+$ ($x\in int(V_+)$), we write $x\geq 0$ (respectively, $x>0$).
For $x\in V$, $x^+$ denotes the projection of $x$ onto $V_+$, and  we let
$x^{-}:=x^+-x$, $|x|:=x^{+}+x^-$. These can also be described via the spectral decomposition  $x=\sum_{1}^{n}x_ie_i$ (where $x_1,x_2,\ldots, x_n$ are the eigenvalues of $x$ and $\{e_1,e_2,\ldots, e_n\}$ is  Jordan frame): $x^+=\sum_{1}^{n}x_{i}^+e_i$,
$|x|=\sum_{1}^{n}|x_i|e_i$, etc.
We see that $|x|^2=x^2$, $\sqrt{x^2}=|x|$,
$\langle x^+,x^-\rangle =0$ and $x^+\circ x^-=0.$
For $x,y\in V$, we define
$$x\sqcap y:=x-(x-y)^+.$$
When $V=\mathbb{R}^{n}$ (with the usual componentwise product and the inner product), this reduces to $\min\{x,y\}$, the componentwise minimum of (vectors) $x$ and $y$ in $\mathbb{R}^{n}$.
For this reason, we may call  the map $(x,y)\rightarrow x\sqcap y$, the `min map' on $V$. The map
$$(x,y)\rightarrow x+y-\sqrt{x^2+y^2}$$ is called the Fischer-Burmeister map. It has been extensively used in the complementarity literature. Below, we state some basic properties of these two maps.

\begin{proposition} \label{basic prop} The following statements hold in $V$:
\begin{itemize}
\item [$(i)$] $u+x\sqcap y=(u+x)\sqcap (u+y).$
\item [$(ii)$] $\lambda(x\sqcap y)=\lambda x\sqcap\lambda y$ for all $\lambda\geq 0$.
\item [$(iii)$] The following are equivalent:
\begin{itemize}
\item [$(a)$] $x\sqcap y=0.$
\item [$(b)$]  $x\geq 0,\,y\geq 0,\,\,\,\mbox{and}\,\,\langle x,y\rangle =0$.
\item [$(c)$] 
$x\geq 0,\,y\geq 0,\,\,\,\mbox{and}\,\, x\circ y=0.$
\end{itemize}
 Moreover, in each case, $x$ and $y$ operator commute.
\item [(iv)] When  $w\geq 0$, the following are equivalent:
\begin{itemize}
\item [$(a)$] $x+y-\sqrt{x^2+y^2+2w}=0$
\item [$(b)$] $x\geq 0,\,y\geq 0,\,\,\,\mbox{and}\,\, x\circ y=w.$ 
\end{itemize}
Moreover, when $w=0$ or $w=e$ (the unit element of $V$), above $x$ and $y$ operator commute.
\end{itemize}
\end{proposition}

\noindent{\bf Proof.} Items $(i)$ and $(ii)$ follow easily from the definition of `min map'. Item $(iii)$ appears in \cite{gowda-sznajder-tao}, Proposition 6 and Item $(iv)$ for $w=0$ or $w=e$ is covered in \cite{gowda-sznajder-tao}, Propositions 6 and 7. Now, let $w\geq 0$ and suppose $x+y-\sqrt{x^2+y^2+2w}=0$. Then, $x+y=\sqrt{x^2+y^2+2w}$. 
This shows that $x+y\geq 0$ and (upon squaring and simplifying) $x\circ y=w$. We need only show that 
$x\geq 0$ and $y\geq 0$. Consider the spectral expansion  $x=\lambda_1e_1+\lambda_2e_2+\cdots+\lambda_ne_n$, where $\lambda_1,\lambda_2,\ldots, \lambda_n$ are the eigenvalues of $x$ and  $\{e_1,e_2,\ldots,e_n\}$ is a Jordan frame in $V$. Suppose, if possible, $x\not\geq 0$; we may assume without loss of generality that $\lambda_1<0$. Then, $x\circ e_1=\lambda_1e_1$ and    

$$0 \leq \langle x+y, e_1 \rangle = \langle x, e_1 \rangle +  \langle y, e_1 \rangle=
\lambda_1 ||e_1||^2+ \frac{1}{\lambda_1}\langle y, x \circ e_1 \rangle  $$
$$= \lambda_1 ||e_1||^2+ \frac{1}{\lambda_1}\langle x \circ y, e_1 \rangle=\lambda_1 ||e_1||^2+\frac{1}{\lambda_1}\langle w, e_1 \rangle<0,$$
as $\langle w,e_1\rangle \geq 0$. This contradiction proves that all eigenvalues of $x$ are nonnegative; so $x\geq 0$.
Similarly, $y\geq 0$. 
Thus we have $(iv)$.
$\hfill$ $\qed$\\

Item $(iv)$ in the above proposition will allow us to formulate a wHLCP as a system of equations. In fact,  $(x,y)$ is a solution of wHLCP$(A,B,w,q)$ (\ref{whlcp}) if and only if it is a solution of the system
\begin{eqnarray*}
\begin{array}{c}
 x+y-\sqrt{x^{2}+y^{2}+2w}=0,\\
 Ax+By-q=0.
\end{array}
\end{eqnarray*}

Our next key result will be used to show that the min  and the Fischer-Burmeister maps are `homotopic'. This will allow us to replace the Fischer-Burmeister map by the `simpler' min map in our main solution analysis.

\begin{proposition}\label{FB prop}
Let $x,y\in V$ and $0\leq t\leq 1$. Then,
$$t\left [x+y-\sqrt{x^2+y^2}\right ]+(1-t)\,x\sqcap y=0\Longleftrightarrow x\sqcap y=0.$$
\end{proposition}

\noindent{\bf Proof.}
In view of Items $(iii)$ and $(iv)$ in the previous proposition, we prove only the `if' part. We also  assume without loss of generality, $0<t<1$. Let $u:=t\left [x+y-\sqrt{x^2+y^2}\right ].$ From Item $(i)$ of the  previous proposition,
$$\Big [(1-t)x+u\Big ]\sqcap \Big [(1-t)y+u\Big ]=0.$$
This implies that $(1-t)x+u\geq 0$ and  $(1-t)y+u\geq 0$.
Now, $(1-t)x+u\geq 0$ implies that $(1-t)x+t\left [x+y-\sqrt{x^2+y^2}\right ]\geq 0$, that is,
$x+ty\geq t\sqrt{x^2+y^2}.$ As $\sqrt{x^2+y^2}\geq \sqrt{y^2}=|y|$ (which is a consequence of the so-called
L\"{o}wner-Heinz inequality, see \cite{gowda-sznajder-tao}, Proposition 8), we see that
$$x\geq t|y|-ty\geq 0.$$
Hence $x\geq 0$. Similarly, $y\geq 0$.
It follows that $\langle x,y\rangle \geq 0$. We now show that $\langle x,y\rangle \leq 0$ to conclude that $\langle x,y\rangle =0$.\\
We first note that
$x\sqcap y = (x+y-|x-y|)/2$.
Let
$$p:=\frac{t}{1-t}\Big [x+y-\sqrt{x^2+y^2}\Big]$$ so that
$x+y-\sqrt{x^2+y^2}=\alpha p$, where $\alpha:=\frac{1-t}{t}$. Then, $\sqrt{x^2+y^2}=(x+y)-\alpha\,p.$ Squaring both sides and simplifying, we get
\begin{equation}\label{intermediate step}
p\circ (x+y)=\frac{1}{2\alpha}\Big [2\,x\circ y+\alpha^2p^2\Big].
\end{equation}
As
$t\left [x+y-\sqrt{x^2+y^2}\right ]+(1-t)\,x\sqcap y=0$, we have $p+x\sqcap y=0,$ that is,
$$2p+(x+y)=|x-y|.$$
Squaring both sides, noting $|x-y|^2=(x-y)^2$, and simplifying, we get
$$4p^2+2 x\circ y+4p\circ (x+y)=-2x\circ y.$$
We replace $4p\circ (x+y)$ by using (\ref{intermediate step}) to get an expression of the form
$$\beta\,x\circ y+\gamma\,p^2=0,$$
where numbers $\beta$ and $\gamma$ are positive. This yields $x\circ y\leq 0$ and $\langle x,y\rangle=\langle x\circ y,e\rangle\leq 0.$
Finally,
since $\langle x, y\rangle\geq 0,$
we have $\langle x, y\rangle= 0.$
Thus we have shown that $x,y\geq 0$ and $\langle x,y\rangle =0.$ Hence, $x\sqcap y=0.$

$\hfill$ $\qed$

We end this subsection by quoting  a well-known determinantal formula.
\begin{proposition}\label{determinant prop}
 \cite{ouelette} For $A,B,X,Y\in\mathbb{R}^{n\times n},$ with $X,Y$ commuting,
the following formula holds:
\begin{eqnarray*}
\det\left[\begin{array}{cc}
A&-B\\
X&Y
\end{array}\right]
=\det(AY+BX).
\end{eqnarray*}
\end{proposition}

A similar statement can be made about linear transformations.
\subsection{Degree theory}
We employ degree theoretic arguments to prove our existence results.
All necessary results concerning degree theory are given in \cite{facchinei-pang}
(specifically Proposition 2.1.3); see also, \cite{lloyd}, \cite{ortega-rheinboldt}.
Here is a brief summary. Suppose $\Omega$ is a bounded open set in $\mathbb{R}^{n}$,
 $g:\overline{\Omega}\rightarrow\mathbb{R}^{n}$
is continuous and $p\not\in g(\partial\Omega)$, where
$\overline{\Omega}$ and $\partial\Omega$
denote, respectively, the closure and boundary of $\Omega.$
Then the degree of $g$ over $\Omega$ with respect to $p$ is defined;
it is an integer and will be denoted by deg$(g,\Omega,p)$. When
this degree is nonzero, the equation $g(x)=p$ has a solution in $\Omega.$
Suppose $g(x)=p$ has a unique solution, say, $x^{*}$ in $\Omega.$
Then deg$(g,\Omega',p)$, which equals deg$(g,\Omega',g(x^*))$, is constant over all bounded open sets $\Omega'$
containing $x^{*}$ and contained in $\Omega.$
This common degree is called the (topological) index of $g$ at $x^{*}$;
 it will be denoted by ind$(g,x^{*})$.
In particular, if $g:\mathbb{R}^{n}\rightarrow\mathbb{R}^{n}$ is a continuous map such
that $g(x) = 0\Leftrightarrow x = 0$, then for any bounded open set containing 0, we have
$$\textrm{ind}(g,0)= \textrm{deg}(g,\Omega,0);$$
moreover, when $g$ is the identity map, $\textrm{ind}(g,0)=1.$

Let $H(x,t):\mathbb{R}^{n}\times[0,1]\rightarrow\mathbb{R}^{n}$ be continuous (in which
case, we say that $H$ is a homotopy). Suppose that for some bounded open set
$\Omega$ in $\mathbb{R}^{n},$ $0\not\in H(\partial\Omega,t)$ for all $t\in[0,1].$¡¡
Then, the \emph{homotopy invariance property of degree} says that
$\textrm{deg}\Big (H(\cdot,t),\Omega,0\Big )$ is independent of $t.$
In particular, if the zero set
\begin{eqnarray*}
\Big \{x: H(x,t)=0\ \textrm{for}\ \textrm{some}\ t\in[0,1]\Big \}
\end{eqnarray*}
is bounded, then for any bounded open set $\Omega$ in $\mathbb{R}^{n}$ containing this zero set,
we have
\begin{eqnarray*}
\textrm{deg}\Big (H(\cdot,1),\Omega,0\Big )
=\textrm{deg}\Big (H(\cdot,0),\Omega,0\Big ).
\end{eqnarray*}

\noindent{\bf Note:} All degree theory concepts and results are also valid over any finite dimensional real Hilbert space (such as $V$ or $V\times V$) instead of $\mathbb{R}^{n}$.
\subsection{A normalization argument}
To show that the zero set of a map or a system of equations is bounded,
we frequently employ the so-called {\it normalization argument.} Here,  a certain sequence of
vectors (with their norms going to infinity) is normalized to yield a unit vector that violates a given criteria.
We illustrate this in the following result, which will be used later.

\begin{proposition}\label{normalization example}
Let $A$ and $B$ two linear transformations on $V$ and $p\in V$. Suppose that
$$\Big [x\sqcap y=0,\,\,Ax+By=0\Big ]\Rightarrow (x,y)=(0,0).$$
Then, the set
$$\Big \{(x,y): x\sqcap y=0,\,\,Ax+By-tp=0\,\,\mbox{for some}\,\,t\in [0,1]\Big \}$$
is bounded.
\end{proposition}

\noindent{\bf Proof.} Suppose the above set is unbounded.
Let $z_k:=(x_k,y_k)$, $t_k\in [0,1]$ with $||z_k||\rightarrow \infty$, and
$x_k\sqcap y_k=0,\,\mbox{and}\,\,Ax_k+By_k-t_kp=0$ for all $k=1,2,\ldots$.
We divide each of the above equations by $||z_k||$ (so as to create normalized vectors $\frac{z_k}{||z_k||}$). We let $k\rightarrow \infty$ and suppose without loss of generality, $x_0:=\lim\frac{x_k}{||z_k||}$ and $y_0:=\lim\frac{y_k}{||z_k||}$. Then, $x_0\sqcap y_0=0,\,\mbox{and}\,\,Ax_0+By_0=0.$ However, from 
$||z_k||^2=||x_k||^2+||y_k||^2$, we get $||x_0||^2+||y_0||^2=1$, which contradicts our assumption. The stated conclusion follows.
$\hfill$ $\qed$

\section{The degree of an ${\bf R}_0$-pair}
In the setting of $V=\mathbb{R}^{n}$, the concepts of LCP-degree of a matrix and HLCP-degree of a pair of matrices 
are useful in describing the existence and stability of solutions, see \cite{gowda-degree} and \cite{sznajder}. 
In what follows, we extend these to Euclidean Jordan algebras. Consider  linear transformations  $M$, $A$, and $B$  on (a general algebra) $V$ and recall that 
$$
\mathrm{HLCP}(A,B,q):=\mathrm{wHLCP}(A,B,0,q)\quad\mbox{and}\quad \mathrm{LCP}(M,q):=\mathrm{HLCP}(I,-M,q).$$
In view of Item $(iii)$ in Proposition \ref{basic prop}, $\mathrm{HLCP}(A,B,q)$ is equivalent to finding $(x,y)\in V\times V$ such that
\begin{eqnarray}
\begin{array}{c}
 x\sqcap y=0,\\
 Ax+By=q,
 \end{array}
\end{eqnarray}
and $\mathrm{LCP}(M,q)$ is equivalent to finding an $x\in V$ such that 
$$x\sqcap (Mx+q)=0.$$
 
We say that $M$ has the ${\bf R}_0$ property on $V$ if
zero is the only solution of LCP$(M,0)$, or equivalently,
$$x\sqcap Mx=0\Leftrightarrow x=0.$$
When this condition holds, for any bounded open set $\Omega$ in $V$ that contains zero,
deg$(\theta,\Omega,0)$ is well defined, where $$\theta(x):=x\sqcap Mx.$$
This common value -- which is ind$(\theta,0)$-- defines the LCP-{\it degree of M}, denoted by $\deg(M)$. 
We now extend this concept to a pair of transformations.

\begin{definition}\label{defn of rnot pair}
The pair $\{A,B\}$ is said to be
an  ${\bf R}_{0}$-pair if
zero is the only solution of HLCP$(A,B,0)$. This means that
$$\Theta(z)=0\Leftrightarrow z=0,$$
where
\begin{eqnarray*}
\Theta(z):=\left[\begin{array}{c}
x\sqcap y\\
Ax+By
\end{array}\right]
\end{eqnarray*}
with $z=(x,y)$.
When this condition holds, ind$(\Theta,0)$ is well defined. (Note that this equals $\deg(\Theta,\Omega,0)$ for any bounded open set $\Omega$ in $V\times V$ that contains zero.)
We define the {\it \rm{HLCP}-degree} of the pair $\{A,B\}$ by
$$\deg(A,B):=\mathrm{ind}(\Theta,0).$$
\end{definition}

Our first result extends Corollary 5.2.6 in \cite{sznajder} from $\mathbb{R}^{n}$ to a general Euclidean Jordan algebra. 

\begin{proposition}
Let $M$ be a linear transformation on $V$ with the ${\bf R}_0$ property.
Then $\{I,-M\}$ is an ${\bf R}_0$-pair and 
$$\deg(I,-M)=\deg(M).$$
\end{proposition}

\noindent{\bf Proof.}
It is easy to see that $\{I,-M\}$ is an ${\bf R}_0$-pair. Now define the map
\begin{eqnarray*}
\Theta(z,t):=\left[\begin{array}{c}
y\sqcap [tx+(1-t)My]\\
x-tMy
\end{array}\right],
\end{eqnarray*}
where $z=(x,y)$ and $t\in [0,1]$. Then,
$$\Theta(z,0)=\left[\begin{array}{c}
\theta(y)\\
x
\end{array}\right]\quad
\mbox{and}\quad \Theta(z,1)=\Theta(z),$$
where  $\theta(y)=y\,\sqcap\, My.$
Moreover, since $M$ has the  ${\bf R}_0$ property (so that $y\sqcap\, My=0\Rightarrow y=0$), it is easy to verify that
$\Theta(z,t)=(0,0)\in V\times V$ if and only if $z=0$. This means that the zero set of $\Theta(z,t)$ (as $t$ varies over $[0,1]$) is just
$\{(0,0)\}$. Thus, for arbitrary  bounded open sets $\Omega_1$ and $\Omega_2$ both containing zero in $V$,
letting $\Omega=\Omega_1\times \Omega_2$, we have,  by the homotopy invariance of degree and the Cartesian product formula (see \cite{facchinei-pang}, Proposition 2.1.3(h)),
\begin{eqnarray*}
\begin{array}{rl}
\deg(I,-M)&=\deg\Big (\Theta(\cdot,1),\Omega,0\Big )=\deg\Big (\Theta(\cdot,0),\Omega,0\Big )\\
&=\deg\Big (\theta,\Omega_1,0\Big )\,\deg\Big (I,\Omega_2,0\Big )=\mathrm{ind}(\theta,0)=\deg(M),
\end{array}
\end{eqnarray*}
where $I$ denotes the identity transformation.
$\hfill$ $\qed$
\\

In the standard LCP theory, ${\bf R}$-matrices \cite{cottle et al} form an important subclass of matrices for which LCP-degree is nonzero. (We note that there are other matrices, such as ${\bf N}$-matrices of first category satisfying this property \cite{gowda-degree}.)
Recall that $M$ is an ${\bf R}$-matrix (in the standard LCP setting) if there is some $d>0$ in $\mathbb{R}^{n}$ such that
zero is the only vector that solves the problems
LCP$(M,0)$ and LCP$(M,d)$. We now consider a generalization.

\begin{definition}\label{defn of R-pair}
Let $A$ and $B$ be two linear transformations on $V$. We say that $\{A,B\}$ is an ${\bf R}$-pair if it  is an ${\bf R}_{0}$-pair and
there exists $p\in V$ such that
\begin{itemize}
\item [$(a)$]  HLCP$(A,B,p)$ has a unique solution, say,  $(\overline{x},\overline{y})$,
\item [$(b)$] $\overline{x}+\overline{y}>0$, and
\item [$(c)$] The derivative of $G$ at $(\overline{x},\overline{y})$ is nonsingular, where
\begin{eqnarray*}
G(z):=\left[\begin{array}{c}
x\sqcap y\\
Ax+By-p
\end{array}\right]
\end{eqnarray*}
with $z=(x,y)$.
\end{itemize}
\end{definition}

Note that condition $(c)$ above is equivalent to: The derivative of $\Theta$ (as given in the Definition of ${\bf R}_0$-pair) at $(\overline{x},\overline{y})$ is nonsingular.
\\

We elaborate on the ${\bf R}$-pair property and give some examples.
Suppose $\{A,B\}$ is an ${\bf R}$-pair.
As $\overline{x}\sqcap \overline{y}=0$ and $\overline{x}+\overline{y}>0$, $\overline{x}$ and $\overline{y}$ operator commute (see Proposition \ref{basic prop}) and so, with respect to some Jordan frame $\{e_1,e_2,\ldots, e_n\}$, we can write
$$\overline{x}=\sum_{1}^{k}\overline{x}_ie_i\quad\mbox{and}\quad \overline{y}=\sum_{k+1}^{n}\overline{y}_je_j,$$
where $1\leq k\leq n$ and $\overline{x}_i,\overline{y}_j>0$ for all $i$ and $j$. Let
$\overline{z}:=(\overline{x},\overline{y})$,
$$\alpha:=\{1,2,\ldots,k\}\quad\mbox{and}\quad \beta=\{k+1,\ldots, n\}.$$
(Note that one of these sets may be empty.)
Thus, the element 
$$\overline{x}-\overline{y}=\sum_{1}^{k}\overline{x}_ie_i-\sum_{k+1}^{n}\overline{y}_ie_i$$
 is invertible (which means that all the eigenvalues are nonzero). In view of Lemma 19 in \cite{gowda-sznajder-tao}, the map $G$ defined in condition $(c)$ is Fr\'{e}chet differentiable.
Let
$$\Phi(x,y):=x\sqcap y=x-\Pi_{V_+}(x-y),$$ where $\Pi_{V_+}$ denotes the projection operator onto $V_+$. Then, the partial derivative of $\Phi$ with respect to $x$ at
$\overline{z}=(\overline{x},\overline{y})$ is given by
$$\Phi_{x}^\prime(\overline{z})=I_x-\Pi_{V_+}^\prime(\overline{x}-\overline{y})\circ I_x,$$
where $I_x$ denotes the identity transformation and $`\circ'$ denotes the composition. Now, we use the formula for the derivative of $\Pi_{V_+}$ given in Lemma 19 of \cite{gowda-sznajder-tao}. (Although this formula is stated in the setting of a  simple Euclidean Jordan algebra, by  writing a general Euclidean Jordan algebra as a product of simple ones, we can show that that the formula is valid 
in any Euclidean Jordan algebra.) Then,  for any $h\in V$ with Peirce decomposition $h=\sum_{1}^{n}h_ie_i+\sum_{i<j}h_{ij}$,
we have
$$\Phi_{x}^\prime(\overline{z})h=h-\Big ( \sum_{1}^{k}h_ie_i+\sum_{i,j\in \alpha;\, i<j}h_{ij}+
\sum_{i\in \alpha,j\in \beta}\frac{\overline{x}_i}{\overline{x}_i+\overline{y}_j}h_{ij}\Big ).$$
Simplification leads to
$$\Phi_{x}^\prime(\overline{z})h=\sum_{k+1}^{n}h_ie_i+\sum_{i,j\in \beta;\, i<j}h_{ij}+\sum_{i\in \alpha, j\in \beta}\frac{\overline{y}_j}{\overline{x}_i+\overline{y}_j}h_{ij}.$$
Similarly, using $\Phi(x,y):=y\sqcap x =y-\Pi_{V_+}(y-x)$, we see that the $y$-derivative of $\Phi(x,y)$ at $(\overline{x},\overline{y})$ is given by
$$\Phi_{y}^\prime(\overline{z})h=\sum_{1}^{k}h_ie_i+\sum_{i,j\in \alpha;\,i<j}h_{ij}+\sum_{i\in \alpha, j\in \beta}\frac{\overline{x}_i}{\overline{x}_i+\overline{y}_j}h_{ij}.$$
Thus, the nonsingularity of $G$ at $(\overline{x},\overline{y})$ could be expressed by  the implication
\begin{equation}
\Big [ \Phi_{x}^\prime(\overline{z})u+\Phi_{y}^\prime(\overline{z})v=0\,\,\mbox{and}\,\,Au+Bv=0\Big]\Rightarrow (u,v)=(0,0).
\end{equation}
Upon simplification, the above condition could be written as
\begin{equation}\label{condition for u and v}
\left.
\begin{array}{c}
u_i=0\,(i\in \beta),\,\,v_j=0\,(j\in \alpha)\\
u_{ij}=0\,(i,j\in \beta;\,i<j),\,\,v_{ij}=0\,(i,j\in \alpha;\,i<j)\\
u_{ij}=-\frac{\overline{x}_i}{\overline{y}_j}v_{ij}\,\,(i\in \alpha,j\in \beta)\\
Au+Bv=0
\end{array}\right\}\Rightarrow (u,v)=(0,0).
\end{equation}

We now show that in some standard settings, condition $(c)$ in Definition \ref{defn of R-pair} 
is a consequence of conditions $(a)$ and $(b)$.
\\

\noindent{\bf Example 1} Let $V=\mathbb{R}^{n}$. Then the standard coordinate vectors form the only Jordan frame and for any element $h\in \mathbb{R}^{n}$, $h_{ij}=0$ for all $i<j$. Thus, in order to verify the implication (\ref{condition for u and v}), we let
$$u=(u_1,u_2,\ldots, u_k,0,0,\ldots,0)^T\quad\mbox{and}\quad
v=(0,0,\ldots,0, v_{k+1},\ldots, v_n)^T $$ and suppose that $(u,v)\neq (0,0)$.
As $\overline{x}=(\overline{x}_1,\overline{x}_2,\ldots,\overline{x}_k,0,0,\ldots,0)^T$ and $\overline{y}=(0,0,\ldots, 0,\overline{y}_{k+1},\ldots \overline{y}_n)^T$, where $\overline{x}_i$ and $\overline{y}_j$ are positive, we see that
for all small $\varepsilon>0$, the pair $(\overline{x}+\varepsilon\,u,\overline{y}+\varepsilon\,v)$ is a solution of HLCP$(A,B,p)$. This contradicts the uniqueness assumption $(a)$. Hence  condition $(c)$ is superfluous in this setting.
\\

\noindent{\bf Example 2} Let $M:V\rightarrow V$ be a linear transformation  that has the
${\bf R}$ property with respect to $V_+$. This means that for some $d>0$ in $V$, zero is the only solution of  the linear complementarity problems
LCP$(M,0)$ and LCP$(M,d)$.
We claim that $\{I,-M\}$ is an ${\bf R}$-pair. It is easy to see that
 the problem HLCP$(I,-M,0)$ has $(0,0)$ as the only solution which means that $\{I,-M\}$ is an ${\bf R}_0$-pair. Also,
 HLCP$(I,-M,d)$ has $(d,0)$ as the only solution. This means that with $(\overline{x},\overline{y})=(d,0)$,
conditions $(a)$ and $(b)$ in the above definition are satisfied. We show that condition $(c)$ holds. If $z=(x,y)$ is close to $(d,0)$, then $x-y$ is close to $d-0$; hence for all such $(x,y)$, $\Pi_{V_+}(x-y)=x-y$ and $x\sqcap y=x-(x-y)=y$. Thus, when $z$ is close to $(d,0)$, we have
$$
G(z):=\left[\begin{array}{c}
 y\\
Ix-My-p
\end{array}\right]
\quad\mbox{and}\quad
G^{\prime}(\overline{z}):=\left[\begin{array}{cc}
0 & I \\
I & -M
\end{array}\right].$$
In view of Proposition \ref{determinant prop}, $G^{\prime}(\overline{z})$ is nonsingular.
Thus, we have verified condition $(c)$. Hence $\{I,-M\}$ is an ${\bf R}$-pair.

\begin{proposition}\label{R implies nonzero degree}
Suppose $\{A,B\}$ is an ${\bf R}$-pair. Then, $\deg(A,B)$ is nonzero.
\end{proposition}

\noindent{\bf Proof.} The pair $\{A,B\}$ satisfies conditions in Definition \ref{defn of R-pair}. Let
\begin{eqnarray*}
G(z,t):=\left[\begin{array}{c}
x\sqcap y\\
Ax+By-tp
\end{array}\right],
\end{eqnarray*}
where $z=(x,y)$ and $t\in [0,1]$.
Then,
$$G(z,1)=\left[\begin{array}{c}
x\sqcap y\\
Ax+By-p
\end{array}\right]\quad
\mbox{and}\quad
G(z,0)=\left[\begin{array}{c}
x\sqcap y\\
Ax+By
\end{array}\right].
$$
As $\{A,B\}$ is an ${\bf R}_0$-pair, by   a normalization argument (see  Proposition \ref{normalization example}), we see that the zero sets of $G(z,t)$ as $t$ varies are uniformly bounded. Suppose $\Omega$  is a bounded open set in $V\times V$ that contains all these zero sets. Note that $G(z,1)$ vanishes only at $(\overline{x},\overline{y})\in \Omega$ and its derivative at this point is nonsingular. Thus,
$$\deg(A,B)=\deg \Big (G(\cdot,0),\Omega,0\Big )=\deg \Big ( G(\cdot,1),\Omega,0\Big )=\mathrm{sgn}\,\det G^{\prime}(\overline{z},1)\neq 0.$$
This proves that $\deg(A,B)$ is nonzero.
$\hfill$ $\qed$

\section{The main existence result}

We now discuss the solvability of wHLCP$(A,B,w,q)$. We recall that $(x,y)$ is a solution of wHLCP$(A,B,w,q)$ if and only if it is a solution of the system
\begin{eqnarray*}
\begin{array}{c}
 x+y-\sqrt{x^{2}+y^{2}+2w}=0,\\
 Ax+By-q=0.
\end{array}
\end{eqnarray*}
We show that this  system has a solution under a nonzero degree condition.

\begin{theorem}\label{main theorem}
Let $\{A,B\}$ be an ${\bf R}_0$-pair with $\deg(A,B)$ nonzero. Then for any  $(w,q)\in V_+\times V$,
$\mathrm{wHLCP}(A,B,,w,q)$ has a nonempty compact solution set.
\end{theorem}

\noindent{\bf Proof}  We fix $(w,q)\in V_+\times V$.
With $z=(x,y)\in V\times V$ and $t\in [0,1]$, we define the following maps:
\begin{eqnarray*}
F(z,t):=\left[\begin{array}{c}
 x+y-\sqrt{x^{2}+y^{2}+2tw}\\
 Ax+By-tq
\end{array}\right],
\end{eqnarray*}
\begin{eqnarray*}
H(z,t):=\left[\begin{array}{c}
t\left [x+y-\sqrt{x^2+y^2}\right ]+(1-t)x\sqcap y\\
 Ax+By
\end{array}\right].
\end{eqnarray*}
We show below that there is some bounded open set $\Omega$ in $V\times V$ which contains all the zeros (in $z$) of
$F$ and $H$ (as $t$ varies over $[0,1]$). Then, over $\Omega$,  $F$ is a homotopy connecting
$F(z,1)$ and $F(z,0)$; $H$ is a homotopy connecting
$H(z,1)$ $\Big (=F(z,0)\Big )$ and $H(z,0)$.
Using the homotopy invariance property of degree, we see that
\begin{eqnarray*}
\textrm{deg}\Big (F(\cdot,1),\Omega,0\Big )
=\textrm{deg}\Big (H(\cdot,0),\Omega,0\Big )\neq 0.
\end{eqnarray*}
This shows that the equation $F(z,1)=0$ has a nonempty bounded solution set.

To justify these, we proceed as follows.
\\
Let
$$Z:=\{z: F(z,t)=0\,\,\mbox{for some}\,\,t\in [0,1]\}.$$
We show by a {\it normalization argument} that $Z$ is bounded. Suppose, if possible, $Z$ is unbounded. Let $z_k:=(x_k,y_k)$, $t_k\in [0,1]$ with $||z_k||\rightarrow \infty$, and
$F(z_k,t_k)=0$ for all $k=1,2,\ldots$.
Let $k\rightarrow \infty$ and without loss of generality, $x_0:=\lim\frac{x_k}{||z_k||}$ and $y_0:=\lim\frac{y_k}{||z_k||}$. We note that
$||x_0||^2+||y_0||^2=1.$
Dividing each component of $F(z_k,t_k)$ by $||z_k||$ and letting $k\rightarrow \infty$, we get
\begin{eqnarray*}
\begin{array}{c}
 x_0+y_0-\sqrt{x_0^{2}+y_0^{2}}=0,\\
 Ax_0+By_0=0.
\end{array}
\end{eqnarray*}
By Proposition \ref{basic prop}, $(x_0,y_0)$ becomes a nonzero solution of HLCP$(A,B,0)$ contradicting the ${\bf R}_0$ property of $\{A,B\}$. Hence, $Z$ is bounded.\\
Next, in  view of Proposition \ref{FB prop} and the ${\bf R}_0$ property of  $\{A,B\}$,
$$\{z: H(z,t)=0\,\,\mbox{for some}\,\,t\in [0,1]\}=\{(0,0)\}.$$
Let $\Omega$ be a bounded open set in $V\times V$ that contains the zero sets of $F$ and  $H$. Then, by the homotopy invariance property of the degree,
$$\textrm{deg}\Big (F(\cdot, 1),\Omega,0\Big )=\textrm{deg}\Big (F(\cdot, 0),\Omega,0\Big )=\textrm{deg}\Big (H(\cdot, 1),\Omega,0\Big )=\textrm{deg}\Big (H(\cdot, 0),\Omega,0\Big )$$ and
$$\textrm{deg}\Big (H(\cdot, 0),\Omega,0\Big )=\textrm{deg}\Big (\Theta,\Omega,0\Big )=\mathrm{ind}(\Theta,0)=\deg(A,B).$$
As the last quantity, by assumption, is nonzero, we conclude that
$$\textrm{deg}\Big (F(\cdot, 1),\Omega,0\Big )\neq 0.$$
This means that the equation $F(z, 1)=0$ has a zero in $\Omega$ proving the existence of a solution of wHLCP$(A,B,w,q)$. As all zeros of $F(\cdot, 1)$ are in the bounded set
$\Omega$ and the solution set of wHLCP$(A,B,w,q)$ is clearly closed, we see nonemptyness and compactness of this solution set. This completes the proof.
$\hfill$ $\qed$
\\

Motivated by interior point methods, we consider the case $w>0$. 
First, we make a simple observation:
$$\Big [x\geq 0,\,y\geq 0,\,\,\mbox{and}\,\,x\circ y>0\Big ]\Longrightarrow x>0\,\,\mbox{and}\,\,y>0.$$
This follows from Item $(iv)$, Lemma 2.6 in \cite{yoshise}. Here is a short/different proof. Let 
$x\geq 0$, $y\geq 0$ and $x\circ y=w>0$. Suppose $x\not>0$ so that
zero is an eigenvalue of $x$. This means that there is a primitive idempotent $e_1$ (which belongs to the Jordan frame that appears in the spectral decomposition of $x$) such that $x\circ e_1=0$. But then, $0<\langle w,e_1\rangle=\langle x\circ y,e_1\rangle =\langle y,x\circ e_1\rangle =0$, is a contradiction. Hence, $x>0$ and, similarly, $y>0$.\

\begin{corollary}\label{ips for whlcp}
Let $\{A,B\}$ be an ${\bf R}_0$-pair with $\deg(A,B)$ nonzero. Suppose $w>0$. Then for any $q\in V$, the following `interior point system' has a nonempty compact solution set:
$$x>0,\,\,y>0,\,\, x\circ y =w,\,\,\mbox{and}\,\,Ax+By=q.$$
\end{corollary}

We now specialize the above two results for a single linear transformation.

\begin{corollary}
Let  $M$ be a linear transformation on $V$. Suppose that $M$ has the ${\bf R}_0$ property and $\deg(M)$ is nonzero. Then, for
 all $(w,q)\in V_+\times V$, the weighted linear complementarity problem 
$$x\geq 0,\,\,y\geq 0,\,\,x\circ y=w,\,\,\mbox{and}\,\,y=Mx+q$$
has a nonempty compact solution set. 
In particular, when $w>0$, for any $q\in V$, the `interior point system'
$$x>0,\,\,y>0,\,\, x\circ y =w,\,\,\mbox{and}\,\,y=Mx+q$$
has a nonempty compact solution set.
\end{corollary}

\noindent{\bf Remarks.} In the standard LCP literature, the solvability and uniqueness issues of interior point systems 
are usually addressed for special types of ${\bf P}_0$-matrices (e.g., ${\bf P}_*$-matrices), see  \cite{kojima book},
Lemma 4.3 and Theorem 4.4. In this regard, the above result appears to be new even in the case of $V=\mathbb{R}^{n}$, as it  holds for numerous types of non ${\bf P}_0$-matrices
such as strictly copositive matrices, ${\bf R}$-matrices, and ${\bf N}$-matrices of first category \cite{cottle et al}. 
\section{The solution set behavior}
Fixing the pair $\{A,B\}$, we let
$\textrm{SOL}(w,q)$
denote the solution set of \\wHLCP$(A,B,w,q)$. The following result describes the behavior of the map $(w,q)\mapsto \textrm{SOL}(w,q)$.

\begin{theorem}
Suppose $\{A,B\}$ is an ${\bf R}_0$-pair with $\deg(A,B)$ nonzero. Then, the following statements hold:
\begin{itemize}
\item [$(a)$] The solution map $(w,q)\mapsto \mathrm{SOL}(w,q)$ from $V_+\times V$ to $V_+\times V_+$ is upper semicontinuous.
\item [$(b)$] Let $w_k\geq 0$ for all $k=1,2,\ldots,$ and $w_k\rightarrow w$. Suppose
$(x_k,y_k)\in \mathrm{SOL}(w_k,q)$ for all $k$. Then, the sequence $\{(x_k,y_k)\}$ is bounded and any accumulation point of this sequence solves $\mathrm{wHLCP}(A,B,w,q)$.
\item [$(c)$] Let $w>0$, $t_k\downarrow 0$, and $(x_k,y_k)\in \mathrm{SOL}(t_k\,w,q)$ for all $k$. Then, $x_k>0$ and $y_k>0$ for all $k$, the sequence $\{(x_k,y_k)\}$ is bounded, and any accumulation point of this sequence solves $\mathrm{HLCP}(A,B,q)$.
\end{itemize}
\end{theorem}

\noindent{\bf Proof.} $(a)$ We fix $(w^*,q^*)\in V_+\times V$ and let $\Omega$ be any open set in $V\times V$ containing SOL$(w^*,q^*)$. We show that for all $(w,q)$ near $(w^*,q^*)$, SOL$(w,q)$ is contained in $\Omega$. Assuming the contrary, suppose there is a sequence $\{(w_k,q_k)\}$ converging to $(w^*,q^*)$ such that some solution $(x_k,y_k)$ in $\textrm{SOL}(w_k,q_k)$ belongs to $\Omega^c$ (the complement of $\Omega$). The sequence $\{(x_k,y_k)\}$ has to be bounded; else, a normalization argument (such as the one used in the previous theorem) produces a nonzero solution of
HLCP$(A,B,0)$ contradicting the ${\bf R}_0$ property of the pair $\{A,B\}$. Now, a subsequential limit of the sequence belongs to SOL$(w^*,q^*)$ and at the same time is in $\Omega^c$ (as this set is closed). This contradiction proves the upper semicontinuity property of the solution set.\\
$(b)$ Under the stated assumptions, $(w_k,q)\rightarrow (w,q)$. A normalization argument shows that the sequence $\{(x_k,y_k)\}$ is bounded. Any subsequential limit of this sequence, clearly, belongs to $\textrm{SOL}(w,q)$,
that is, solves wHLCP$(A,B,w,q)$.\\
$(c)$ That $x_k>0$ and $y_k>0$ for all $k$ follows from 
Corollary \ref{ips for whlcp}.
For the remaining statements,  we specialize $(b)$ with $w_k:=t_kw$.
$\hfill$ $\qed$

\section{P-pairs over $V$}
A linear transformation $M$ on $V$ is said to be a ${\bf P}$-transformation \cite{gowda-sznajder-tao} if
\[
\left.\begin{array}{r}
 x\,\,\mbox{and}\,\,Mx\,\,\mbox{operator commute}\\
x\circ Mx\leq0
\end{array}\right\}
\Rightarrow x=0.
\]
${\bf P}$-transformations are generalizations of
${\bf P}$-matrices. An important example of ${\bf P}$-transformations appears in dynamical systems: the Lyapunov transformation $X\mapsto AX+XA^T$ on the Euclidean Jordan algebra
of $n\times n$ real symmetric matrices is a ${\bf P}$-transformation if and only if the (real square) $A$ is positive stable. See \cite{gowda-sznajder-tao} for properties of ${\bf P}$-transformations and further examples. 
We now extend this notion to a pair of transformations.

\begin{definition}
A pair of linear transformations  $\{A,B\}$ is said to be a  ${\bf P}$-pair over $V$ if
\[
\left.\begin{array}{r}
 x\,\,\mbox{and}\,\,y\,\,\mbox{operator commute}\\
x\circ y\leq0\\
 Ax+By=0
\end{array}\right\}
\Rightarrow (x,y)=(0,0).
\]
\end{definition}

Below, we collect some properties of such pairs.

\begin{proposition}
Suppose $\{A,B\}$ is a ${\bf P}$-pair. Then, the following statements hold.
\begin{itemize}
\item [$(a)$] $A$ and $B$ are invertible.
\item [$(b)$] $-B^{-1}A$ and $-A^{-1}B$ are  ${\bf P}$-transformations.
\item [$(c)$]  $\{A,B\}$ is an ${\bf R}$-pair.
\item [$(d)$] For all $(w,q)\in V_+\times V$, $\mathrm{wHLCP}(A,B,w,q)$ has a nonempty compact solution set.
\end{itemize}
\end{proposition}

\noindent{\bf Proof.}
$(a)$ If $Ax=0$ for some $x$, then $$\Big [\, x\,\,\mbox{and}\,\,0\,\,\mbox{operator commute},\,\, x\circ 0=0,\,Ax+B0=0\Big ]\Rightarrow (x,0)=(0,0).$$ This shows that $A$  is invertible. Similarly $B$ is invertible.\\
$(b)$ Let $M:=-B^{-1}A$.
If $x$ and $Mx$ operator commute and $x\circ Mx\leq 0$, then, with $y:=Mx=-B^{-1}Ax$, we see that:
$x$ and $y$ operator commute, $x\circ y\leq 0$, and $Ax+By=0$. Hence $(x,y)=(0,0)$ and so $x=0.$ Thus, $-B^{-1}A$ is a ${\bf P}$-transformation. Similarly, 
$-A^{-1}B$ is also a ${\bf P}$-transformation.\\
$(c)$
We now show that $\{A,B\}$ is an ${\bf R}$-pair.  First, $\{A,B\}$ is an ${\bf R}_0$-pair: When
$x\sqcap y=0$ and $Ax+By=0$, by Proposition 2, $x\,\,\mbox{and}\,\,y\,\,\mbox{operator commute},$
 $x\circ y=0$, and $Ax+By=0$; so,
$(x,y)=(0,0)$ by the definition of a {\bf P}-pair.

Now let $p:=Be$, where $e$ denotes the unit element in $V$. Clearly, $(0,e)$ is a solution of  HLCP$(A,B,p)$. We show that this is the only solution.
Let $(x,y)$ be any solution of HLCP$(A,B,p)$ so that
$x\,\,\mbox{and}\,\,y\,\,\mbox{operator commute}$, $x,y\geq 0$, and $x\circ y=0$. This implies that $x\,\,\mbox{and}\,\,y-e\,\,\mbox{operator commute}\,\,$ and $x\circ (y-e)=-x\leq 0$. Since we also have $Ax+B(y-e)=0$, by the definition of  ${\bf P}$-pair, $(x,y-e)=(0,0)$. This implies that $(x,y)=(0,e)$.
Since $0+e=e>0$, we see that conditions $(a)$ and $(b)$ in Definition \ref{defn of R-pair} 
hold. We now verify condition
$(c)$ in that definition. If $(x,y)$ is close to $(0,e)$, then $x-y$ is close to $-e$ and so $\Pi_{V_+}(x-y)=0.$ In this case, $x\sqcap y=x-\Pi_{V_+}(x-y)=x$.
Hence,
for all $z=(x,y)$ near $\overline{z}:=(0,e)$,
$$
G(z):=\left[\begin{array}{c}
 x\\
Ax+By-Be
\end{array}\right]
\quad\mbox{and}\quad
G^{\prime}(\overline{z}):=\left[\begin{array}{cc}
I & 0 \\
A & B
\end{array}\right].$$
In view of Proposition \ref{determinant prop} and the invertibility of $B$, $G^{\prime}(\overline{z})$ is nonsingular.
\\
$(d)$ This follows from Proposition \ref{R implies nonzero degree}  and Theorem \ref{main theorem}.
$\hfill$ $\qed$
\\

\noindent{\bf Example 3} Let $V={\cal S}^n$, the Euclidean Jordan algebra  of all $n\times n$ real symmetric matrices with 
$\langle X,Y\rangle:=trace(XY)$ and $X\circ Y:=XY+YX$. Here, ${\cal S}^n_+$ is the `semidefinite cone' consisting of positive semidefinite matrices in ${\cal S}^n$. We write 
$X\succeq 0$ and $X\succ 0$, respectively, to denote elements of $V_+$ and its interior. Let $A$ be an $n\times n$ real matrix which is
positive stable (which means that every eigenvalue of $A$ has positive real part). Then, the Lyapunov transformation 
$X\mapsto AX+XA^T$ is a ${\bf P}$-transformation on ${\cal S}^n$ \cite{gowda-sznajder-tao}. Consequently, for any 
$W\succ 0$, the following system has a solution:
$$X\succ 0,\,\,Y\succ 0,\,\, X\circ Y =W,\,\,\mbox{and}\,\,Y=AX+XA^T.$$
While the existence of an $X\succ 0$ with $AX+XA^T\succ 0$ is already covered in the Lyapunov theory of dynamical systems, 
what is new here is that such an $X$ can be found satisfying an additional condition $X\circ (AX+XA^T)=W$ where $W\succ 0$ is arbitrary! 
One can make statements similar to the above for the Stein transformation $X\mapsto X-BXB^T$, where $B$ is an $n\times n$ real  Schur stable matrix
(which means that every  eigenvalue of $B$ has absolute value  less than one) \cite{gowda-sznajder-tao}.
\section{${\bf P}$-pairs over $\mathbb{R}^{n}$}

We now consider $V=\mathbb{R}^{n}$ and prove a uniqueness result.
\begin{theorem}
Let $V=\mathbb{R}^{n}$. Then, the following statements are equivalent:
\begin{itemize}
\item [$(a)$] $\{A,B\}$ is a ${\bf P}$-pair.
\item [$(b)$] $\mathrm{wHLCP}(A,B,w,q)$ has a unique solution  for every  $(w,q)\in\mathbb{R}^{n}_{+}
\times \mathbb{R}^{n}$.
\item [(c)] $\mathrm{HLCP}(A,B,q)$ has a unique solution for every $q\in \mathbb{R}^{n}.$
\end{itemize}
\end{theorem}

\noindent{\bf Proof} $(a)\Rightarrow (b)$: The solvability of wHLCP$(A,B,w,q)$ has been addressed in the previous result.
We now prove uniqueness.
 Suppose that $(x_{1},y_{1})$ and $(x_{2},y_{2})$ are any two solutions of  wHLCP$(A,B,w,q)$,
i.e.,
\begin{eqnarray*}
\begin{array}{rrl}
\left\{\begin{array}{l}
 x_{1}\geq0,\ y_{1}\geq0\\
 x_{1}\ast y_{1}=w\\
 Ax_{1}+By_{1}=q
\end{array}\right.
&\quad\textrm{and}\quad
&\left\{\begin{array}{l}
 x_{2}\geq0,\ y_{2}\geq0\\
 x_{2}\ast y_{2}=w\\
 Ax_{2}+By_{2}=q.
\end{array}\right.
\end{array}
\end{eqnarray*}
As $w\geq 0$, let 
$$\alpha:=\{i:w_i>0\}\quad\mbox{and}\quad \beta:=\{i:w_i=0\}.$$
Then, 
\begin{eqnarray*}
\begin{array}{c}
(x_{1}\ast y_{1})_{i}=(x_{2}\ast y_{2})_{i}=w_{i}>0,\quad \forall\,i\in \alpha,\\
(x_{1}\ast y_{1})_{i}=(x_{2}\ast y_{2})_{i}=w_{i}=0,\quad \forall\,i\in \beta.
\end{array}
\end{eqnarray*}
Now considering the componentwise product (only) over the $\alpha$ indices, we have 
\begin{eqnarray*}
\begin{array}{l}
(x_{1}-x_{2})\ast(y_{1}-y_{2})\\
=x_{1}\ast y_{1}+x_{2}\ast y_{2}-x_{1}\ast y_{2}-x_{2}\ast y_{1}\\
=x_{1}\ast y_{1}+x_{1}\ast y_{1}-x_{1}\ast \left(\dfrac{x_{1}\ast y_{1}}{x_{2}}\right)-x_{2}\ast y_{1}\\
=(x_{1}-x_{2})\ast y_{1}+\left(\dfrac{x_{1}\ast y_{1}}{x_{2}}\right)\ast(x_{2}-x_{1})\\
=(x_{2}-x_{1})\ast\left(\dfrac{x_{1}\ast y_{1}}{x_{2}}-y_{1}\right)\\
=-\dfrac{y_{1}}{x_{2}}(x_{1}-x_{2})^{2}\\
\leq0.
\end{array}
\end{eqnarray*}
And over the $\beta$ indices,
\begin{eqnarray*}
\begin{array}{l}
(x_{1}-x_{2})\ast(y_{1}-y_{2})\\
=x_{1}\ast y_{1}+x_{2}\ast y_{2}-x_{1}\ast y_{2}-x_{2}\ast y_{1}\\
=0+0-x_{1}\ast y_{2}-x_{2}\ast y_{1}\\
\leq0.
\end{array}
\end{eqnarray*}
Therefore, we have
\begin{eqnarray*}
\begin{array}{l}
 (x_{1}-x_{2})\ast (y_{1}-y_{2})\leq0,\\
 A(x_{1}-x_{2})+B(y_{1}-y_{2})=0.
\end{array}
\end{eqnarray*}
As $\{A,B\}$ is a ${\bf P}$-pair and vectors in $\mathbb{R}^{n}$ always operator commute, we see that    $x_{1}=x_{2}$ and $y_{1}=y_{2}.$
Thus we have uniqueness of solution in any wHLCP$(A,B,w,q)$.
\\
$(b)\Rightarrow (c):$ This is obvious by taking $w=0$.\\
$(c)\Rightarrow (a):$
Suppose $(x,y)\neq (0,0)$ with $x\circ y\leq 0$ and $Ax+By=0$. For $q:=Ax^++By^+=Ax^-+By^-$, we see that $(x^+,y^+)$ and $(x^-,y^-)$ are two distinct
solutions of HLCP$(A,B,q)$. This contradicts condition $(c)$.
$\hfill$ $\qed$
\\

We remark that such a uniqueness result may not  prevail over  general Euclidean Jordan algebras even for {\bf P}-transformations, see the remarks following Theorem 14 in \cite{gowda-sznajder-tao}.\\ 
The above result, especially Item $(c)$, allows us to connect ${\bf P}$-pairs to the so-called ${\bf W}$ property for a pair of matrices,
see \cite{sznajder-gowda}. 
\\

\noindent{\bf Concluding Remarks.} In this paper, we have presented some existence and uniqueness results for 
weighted horizontal linear complementarity problems over Euclidean Jordan algebras. These are established for 
${\bf R}_0$-pairs of linear transformations satisfying a (nonzero) degree condition. The novelty here is the use of 
`weighted'  
Fischer-Burmeister map and degree theory techniques. We hope to consider applications, algorithms, and non ${\bf R}_0$-pairs in a future study.

\begin{acknowledgements}
The  work of the first author  is supported by the National Natural Science Foundation of China (No. 11401126)
and Guangxi Natural Science Foundation \\(Nos. 2016GXNSFBA380102, 2014GXNSFFA118001), China. The
third author was supported by Loyola Summer Research Grant 2017.
\end{acknowledgements}



\end{document}